\documentclass[11pt,reqno]{amsart}

\newtheorem{theorem}{Theorem}[section]
\newtheorem{lemma}[theorem]{Lemma}
\newtheorem{corollary}[theorem]{Corollary}
\newtheorem{proposition}[theorem]{Proposition}

\theoremstyle{definition}

\theoremstyle{remark}
\newtheorem{remark}[theorem]{Remark}

\begin{document}

\title[Branching rules for symplectic groups]{Branching rules for              
modular fundamental representations of symplectic groups}

\author{A.~A.~Baranov}
\address{Institute of Mathematics, National Academy of Sciences of Belarus,         
Surganova 11, Minsk, 220072, Belarus.                                        
}
\email{baranov@im.bas-net.by}

\author{I.~D.~Suprunenko}
\address{Institute of Mathematics, National Academy of Sciences of Belarus,         
Surganova 11, Minsk, 220072, Belarus.}
\email{suprunenko@im.bas-net.by}
\thanks{Both authors have been supported 
by 
the Institute of Mathematics of the National                                 
Academy of Sciences of Belarus in the framework                              
of the program ``Mathematical structures''
and by
INTAS (Project 93-183-Ext).}

\date{}

\dedicatory{To our teacher A.\ E.\ Zalesskii on the
occasion of his sixtieth birthday}

\begin{abstract}
In this paper branching rules for the fundamental 
representations of the symplectic groups in positive characteristic are
found. 
The submodule structure of the restrictions of the fundamental modules for 
the group $Sp_{2n}(K)$ to the naturally embedded subgroup $Sp_{2n-2}(K)$ is 
determined. 
As a corollary, inductive systems of
fundamental representations for $Sp_{\infty}(K)$ are classified.
The submodule structure of the fundamental Weyl modules is refined.
\end{abstract}

\newcommand{\cal}{\mathcal}

\newcommand{\Ze}{{\mathbb Z}}
\newcommand{\gP}{\mathfrak P}
\newcommand{\gQ}{\mathfrak Q}

\newcommand{\nmid}{\mathbin{\not|}}
\newcommand{\restr}{{\downarrow}}
\newcommand{\lp}{\mathop{\rm lp}\nolimits}
\newcommand{\dkl}{d_k^{l+2k}}
\newcommand{\om}{\omega}

\maketitle


\section{Introduction}\label{1} 

The article is devoted to finding branching rules for the fundamental representations of the 
sympplectic groups in positive characteristic. The classical branching rules 
are concerned with the restrictions of representations of the classical 
algebraic and symmetric groups 
in characteristic 0 to naturally embedded subgroups of smaller ranks. For a 
group of rank $n$ and its fixed irreducible representation $\varphi$ 
they yield 
the composition factors of the restriction of $\varphi$ to a naturally embedded 
subgroup of rank $n-1$ and hence to similar subgroups of smaller ranks, at 
least algorithmically. These rules provide a basis 
for induction on rank and have found numerous applications. In positive 
characteristic one cannot expect to obtain complete
branching rules in an explicit 
form in a near future since this problem is closely connected 
with that of finding the dimensions of arbitrary irreducible representations
and the composition factors of the Weyl modules.
So it is worth to investigate important particular cases where such rules 
can be found and to seek for asymptotic analogs of these rules. The notion 
of an inductive system of representations (see the definition below)
introduced by Zalesskii in
\cite{Z1}
yields an asymptotic version of the branching rules.
It  proved to be useful for the study of ideals in group algebras of locally
finite groups as well, see, for instance, Zalesskii's survey \cite{Z2}. We 
classify the inductive systems of the fundamental representations for the 
infinite-dimensional symplectic group $Sp_{\infty}(K)$. This class of 
representations yields an example of representations of a simple form for 
which the branching rules in positive characteristic differ from the 
characteristic 0 case. 

      Let $K\subset F$ be fields of characteristic $p>0$, 
$\bar F$ be the algebraic closure of $F$, and
$\Ze^+$ be the set of nonnegative integers. 
Let $G_n=Sp_{2n}(K)$. Denote by 
$\omega_i^n$, $0\le i\le n$, the $i$th fundamental module and representation 
of $G_n$ over $F$ where $\omega_0^n$ is the trivial 
one. Let $W_0^n=\om_0^n, W_1^n,\dots,W_n^n$ be the corresponding Weyl
modules.
Set  $W_i^n=\om_i^n=0$ for $i<0$ and for $i>n$. The labeling
of the fundamental modules is standard, the fundamental and the Weyl modules
for $G_n$ are the $F$-modules 
affording the
restrictions to $G_n$ of the relevant representations of
the group $Sp_{2n}(\bar F)$
(it is well known that these restrictions can be realized over $F$). 
For an integer $z>0$ we denote by $\lp(z)$ the maximal $i$ such
that $p^i\mid z$. We have $\lp(z)=0$ if $p\nmid z$. 
Let $M$ be a $G_n$-module.
The restriction of $M$ to $G_{n-1}$ is denoted by 
$M\restr G_{n-1}$.  
We shall write
$$
M\sim N_1+\dots+N_q
$$ 
if there is a series 
$
0=M_0\subset M_1\subset\dots\subset M_q=M
$
of submodules of $M$ and a permutation $\sigma$ such that
$N_{\sigma(i)}\cong M_i/M_{i-1}$ for all $i=1,\dots,q$. 
Moreover, if in addition
$M_i/M_{i-1}$ coincides with the socle of $M/M_{i-1}$ for $i=1,\dots,q$,
then the sequence
$$
N_{\sigma(1)}\prec_s N_{\sigma(2)}\prec_s \ldots\prec_s N_{\sigma(q)}
$$
is called the {\em socle series} of $M$.

Theorem~\ref{br} below  describes the
branching rules for the fundamental $G_n$-modules and the submodule 
structure of the restrictions of these modules to $G_{n-1}$. 

\begin{theorem}\label{br}
Let 
$n\ge2$ and 
$0\le i \le n$.
Set $d=\lp(n-i+1)$;
$\varepsilon=0$ if $n-i+1\equiv-p^d \pmod{p^{d+1}}$ and 
$\varepsilon=1$ otherwise. 
Then

{\rm(i)}
$\om_i^n\restr G_{n-1}\sim 
\om_i^{n-1}+2\om_{i-1}^{n-1}+\left(\sum_{t=0}^{d-1}2\om_{i-2p^t}^{n-1}\right)+
\varepsilon\om_{i-2p^d}^{n-1}
$
(the sum in the brackets is zero whenever $d=0$);

{\rm(ii)}
$\omega_i^n\restr G_{n-1}=\omega_{i-1}^{n-1}\oplus
\omega_{i-1}^{n-1}\oplus D$
and the series 
$$
\om^{n-1}_{i-2}\prec_s\om^{n-1}_{i-2p}\prec_s\dots
\prec_s\om^{n-1}_{i-2p^{d-1}}\prec_s
\omega
\prec_s\om^{n-1}_{i-2p^{d-1}}\prec_s
\dots\prec_s\om^{n-1}_{i-2p}\prec_s\om^{n-1}_{i-2} 
$$
with $\omega=\om_i^{n-1}\oplus\varepsilon\om_{i-2p^d}^{n-1}$
and $\om_j^{n-1}$ omitted for $j<0$ 
is the socle series of $D$. In particular, 
$D=\omega$ if $i=0,1$ or $p\nmid n+1-i$.
\end{theorem}

\begin{corollary}\label{ccr} For $n\ge2$ the restriction 
$\omega_i^n\restr G_{n-1}$ is completely reducible if and only if 
$i=0,1$ or $p\nmid n+1-i$. 
\end{corollary}

The proof of Theorem~\ref{br} is based on the description of
the composition factors of the fundamental Weyl modules (Premet and
Suprunenko \cite[Theorem~2]{PS} for $p>2$ and independently Adamovich 
\cite[Theorem~2 and its Corollary~1]{A1} for 
arbitrary $p$) and Adamovich's results \cite{A2} 
on the submodule structure 
of these Weyl modules. 
In Section~\ref{2} these results are refined (Theorem~\ref{weylmain}).
In particular, a new irreducibility
criterion for the fundamental Weyl modules
is obtained (Corollary~\ref{irredweyl})
and it is proved that their socles
are always simple (Corollary~\ref{socleweyl}). 

For $n-p+2\le i\le n$ Gow \cite{Go} has given an explicit construction of the modules 
$\omega^n_i$ and has described the submodule 
structure of the restrictions $\omega^n_i\restr G_1\times G_{n-1}$
(the natural embedding) (\cite[Theorem 2.2]{Go}). This implies our 
Theorem~\ref{br} for these modules. In \cite{Go} a certain explicitly 
determined operator $\delta$ on the exterior algebra $\wedge V$ of the 
natural $G_n$-module $V$ is considered and it is proved that for 
$n-p+2\le i\le n$ the module $\omega^n_i$ can be realized as the quotient 
${\rm ker}\,\delta\cap\wedge^iV/\delta^{p-1}(\wedge^{i+2p-2}V)$ 
(\cite[Corollary 2.4]{Go}). This nice construction gives a realization for 
an important class of modules without complicated representation-theoretic 
machinery. However, it cannot be extended to other fundamental modules
since according to \cite[Theorem 4.2]{Go}, the quotient above is zero for 
$i<n-p+2$. 

In \cite{S} Sheth has found the branching rules
for modular representations of symmetric groups corresponding
to two part partitions. The composition factors occurring
in the relevant restrictions are similar to those of the module $D$
in Theorem~\ref{br}(ii). We conjecture that the submodule
structure of these restrictions is also similar to that of $D$.
The authors plan to consider this question as well
as the similar one for representations of special
linear groups 
with highest weights $\omega_i+\omega_j$
in a subsequent paper.

In Section~\ref{4} Theorem~\ref{br} is applied to classify the inductive
systems 
of fundamental $F$-representations for 
$Sp_{\infty}(K)$. Let 
\begin{equation}\label{eqse}
H_1\subset H_2\subset\ldots\subset H_n\subset\ldots
\end{equation}
be a sequence of groups, and $\Psi_n$, $n=1,2, \ldots$, be a nonempty 
finite set of (inequivalent) 
irreducible representations of $H_n$ over a fixed field. The system  
$\Psi=\{\Psi_n\mid n=1,2,\ldots\}$ is called an {\em inductive system} 
(of representations) for the group $H=\bigcup_{n=1}^{\infty} H_n$
if each $\Psi_n$ coincides with the union of the sets of 
composition factors (up to equivalence) 
of the restrictions $\pi\restr H_n$ where $\pi$ runs 
over $\Psi_{n+1}$. 
In this article (\ref{eqse}) is the sequence of the
naturally embedded groups $G_n=Sp_{2n}(K)$,
so $\bigcup_{n=1}^{\infty} G_n=Sp_{\infty}(K)$. Set 
$$
\begin{array}{rclrcl}
{\cal F}_n&=&\{\om^n_i\mid 0\le i\le n\}, 
& {\cal F}&=&\{{\cal F}_n\mid n=1,2, \ldots\}; \\
{\cal L}_n^s&=&\{\om^n_i\mid 0\le i\le s\},  
& {\cal L}^s&=&\{{\cal L}_n^s\mid n=1,2, \ldots\}; \\
{\cal R}_n^u&=&\{\om^n_i\mid n+1-u\le i\le n\},
& {\cal R}^u&=&\{{\cal R}_n^u\mid n=1,2, \ldots\}.
\end{array}
$$

\begin{theorem}\label{ti} The inductive systems of 
fundamental representations over $F$ for $Sp_{\infty}(K)$
are exhausted by the systems  ${\cal F}$, ${\cal L}^s$, ${\cal R}^{p^t-1}$, 
and ${\cal L}^s\cup{\cal R}^{p^t-1}$ ($s\ge0$, $t\ge1$). 
\end{theorem}

It is
clear that ${\cal L}^0$ (which 
consists of the trivial representations) 
and ${\cal R}^{p-1}$ are minimal
inductive systems. However, the question on the minimal inductive systems for
$Sp_{\infty}(K)$ is far from solution. For $p>2$ Zalesskii and
Suprunenko \cite{SZ} have described the inductive  system 
$\Phi=\{\Phi_n\mid n=1,2,\ldots\}$ where for each $n$ the set $\Phi_n$
consists of two irreducible representations with highest weights
$\om_{n-1}+\frac{1}{2}(p-3)\om_n$ and $\frac{1}{2}(p-1)\om_n$. The system 
$\Phi$ coincides with ${\cal R}^2$ for $p=3$ and yields another example of
a minimal inductive system for $p>3$. 

For other classical groups the questions investigated in this paper present  no 
problems since the situation is the same as in characteristic 0.

The authors \cite{BS} have found the minimal and the minimal nontrivial 
inductive systems for the group $SL_{\infty}(\bar K)$. For this group the 
system consisting of the trivial representation is the only minimal inductive 
system, and the minimal nontrivial ones are exhausted by the systems 
${\bf L}^j=\{L^j_n\mid n=1,2,\ldots\}$ and 
${\bf R}^j=\{R^j_n\mid n=1,2,\ldots\}$ where $L^j_n$ consists of two 
irreducible representations of $SL_{n+1}(\bar K)$ with highest weights $0$
and $p^j\omega_1$ and $R^j_n$ of those with highest weights $0$ and 
$p^j\omega_n$. The picture is similar  for the groups $SL_{\infty}$ and 
$SU_{\infty}$ over locally finite fields. 

Until Proposition~\ref{lred} we assume that $K=F=\bar F$.  
At the end Proposition~\ref{lred} transfers the results to arbitrary
fields.


\section{The structure of the 
fundamental Weyl modules}\label{2} 
In this section we refine the results of 
\cite{PS}, \cite{A1}, and \cite{A2} on the structure of the 
fundamental Weyl modules for $G_n$.     

Throughout the paper we set 
$\pi^n_i=\om^n_{n+1-i}$ and $V^n_i=W^n_{n+1-i}$. 
We denote by $[a,b]$ the set of all $j\in\Ze^+$ with $a\le j\le b$.
For an integer $k\in\Ze^+$ write its $p$-adic expansion  
$k=k_0+k_1p+\dots+k_sp^s$ with $0\le k_i<p$ and set $k_i=0$ 
for all such $i\in\Ze^+$ that $p^i>k$.  
We shall write $k=(k_0,k_1,\dots,k_s)$. 
We say that an integer $m$ 
{\em contains $k$ to base $p$}  
and write $k\subset_p m$ 
if and only if for each $i$ either $k_i=0$, or
$k_i=m_i$. Set
$d_k^m=1$ if $k\subset_p m$, and $d_k^m=0$ otherwise.

\begin{theorem}\cite[Theorem~2]{PS}\label{dec} Let $p>2$. Then 
$W_i^n\sim \sum_{k=0}^{\infty}d_k^{n+1-i+2k}\om_{i-2k}^n$.
\end{theorem}

We need some more notation to state Adamovich's results. 
For
$\lambda\in\Ze^+$  define  maps $s_{\lambda}':\Ze^+\to\Ze^+$ and 
$s_{\lambda}:\Ze^+\to\Ze^+$ setting

$s_{\lambda}'(l)=l+2k'$ where
$l+1=a'p^{\lambda}-k'$, $a'\in\Ze^+$, $0\le k'< p^{\lambda}$;
 
$s_{\lambda}(l)=l+2k$ where
$l=ap^{\lambda}-k$, $a\in\Ze^+$, $0\le k< p^{\lambda}$. 

\noindent 
We say that the
{\em reflection} $s_{\lambda}'$ or $s_{\lambda}$ is {\em $l$-admissible} 
if $k'\ne0$ and $p\nmid a'$ or $k\ne0$ and $p\nmid a$, respectively.
We denote by $S(l)$ the set of all $m>l$ that can be
written in the form 
$m=s_{\lambda_u}\dots s_{\lambda_1}(l)$ 
where $\lambda_u<\dots<\lambda_1$ and for each $i=0,1,\dots,u-1$ the
reflection
$s_{\lambda_{i+1}}$ is $s_{\lambda_i}\dots s_{\lambda_1}(l)$-admissible.
Similarly we define $S'(l)$ (writing $s_{\lambda_i}'$
instead of $s_{\lambda_i}$).

\begin{theorem}\cite{A1}\label{Ad1} 
Let $0\le l\le n$. Then 
$V_{l+1}^n\sim\pi_{l+1}^n+\sum_{m\in S'(l)}\pi_{m+1}^n$.
\end{theorem}

As $s_{\lambda}'(x-1)=s_{\lambda}(x)-1$, 
the following theorem yields an equivalent statement. 

\begin{theorem}\label{Ad1n} 
Let $1\le l\le n+1$. Then 
$V_{l}^n\sim\pi_{l}^n+\sum_{m\in S(l)}\pi_{m}^n$.
\end{theorem}

Let us rewrite Theorem~\ref{dec} in terms of $\pi_m^n$
and $V_l^n$ (without restrictions on $p$). 

\begin{theorem}\label{dec2} Let $1\le l\le n+1$. Then 
$V^n_l\sim \sum_{k=0}^{\infty}d_k^{l+2k}\pi^n_{l+2k}$.
\end{theorem}

Now our goal is to show that Theorems~\ref{Ad1n} and \ref{dec2} are
equivalent, so Theorem~\ref{dec2} (and \ref{dec}) holds for $p=2$. 
For this purpose
we prove 
some technical facts on the triples $k, l,m$ with $k\subset_p m=l+2k$ and 
admissible reflections. 

Until the end of the section $l\ge1$. 
 For each $m\in S(l)$ the tuple 
$(\lambda_1;\dots;\lambda_u)$ is uniquely determined (see the comments before
the Theorem in \cite{A2}).  
If $u$ is odd for some $m$, set $\lambda_{u+1}=\lp(m)$.
Then $s_{\lambda_{u+1}}(m)=m$ and $\lambda_{u+1}<\lambda_u$. Now for every 
$m\in S(l)$ we have a uniquely determined sequence of reflections
$s_{\lambda_1}, \ldots, s_{\lambda_{2t}}$. Such sequences will be
called {\em $l$-admissible}. 
For an  integer $0\le a\le p-1$ set $\bar a=p-1-a$. 
The following lemma is straightforward.

\begin{lemma}\label{l28} 
Set $q=\lp(l)$. The reflection $s_{\lambda}$ is $l$-admissible
if and only if  $\lambda>q$ and $l_{\lambda}\ne p-1$.
In that case 
$
s_{\lambda}(l)=p^q(\bar l_q+1,\bar l_{q+1},\dots
\bar l_{\lambda-1},l_{\lambda}+1,l_{\lambda+1},\dots).
$
\end{lemma}

Two consequent applications of Lemma~\ref{l28} yield 

\begin{proposition}\label{l29}
Let $\lp(l)\le\mu<\lambda$,
$m=s_{\mu}s_{\lambda}(l)$, and $k=(m-l)/2$. The pair
$s_{\lambda}, s_{\mu}$ is $l$-admissible if and only if
$l_{\lambda}\ne p-1$ and $l_{\mu}\ne0$. In that case
$$
\begin{array}{rclrrrllllll}
m&= &(l_0,&\dots,&l_{\mu-1},&\bar l_{\mu}+1,&\bar l_{\mu+1},&\dots,&
\bar l_{\lambda-1},&l_{\lambda}+1,&l_{\lambda+1},&\dots), \\
k&= &(0,  &\dots,&        0,&\bar l_{\mu}+1,&\bar l_{\mu+1},&\dots,&
\bar l_{\lambda-1}).
\end{array}
$$
In particular, $k\subset_p m=l+2k$.
\end{proposition}

We call a tuple $\sigma=(\lambda_1;\dots;\lambda_{2t})$
{\em $l$-admissible} if $\lambda_i\in\Ze^+$,
$\lambda_1>\dots>\lambda_{2t}$,
$l_{\lambda_{2j-1}}\ne p-1$ and $l_{\lambda_{2j}}\ne 0$
for $j=1,\dots,t$. 
For an $l$-admissible $\sigma$ set
$$
\gQ_l(\sigma)=\bigcup_{j=1}^t[\lambda_{2j},\lambda_{2j-1}-1],
$$
$\delta_{i\sigma}=1$ if $i=\lambda_j$ for some $j$ and
$\delta_{i\sigma}=0$ otherwise. Define $l^{\sigma}\in\Ze^+$ putting
$$
l^{\sigma}_i=\left\{
\begin{array}{ll}
{\bar l}_i+\delta_{i\sigma}, & i\in\gQ_l(\sigma) \\
l_i+\delta_{i\sigma}, & i\notin\gQ_l(\sigma).
\end{array}\right.
$$ 
Proposition~\ref{l29} yields the following corollary.

\begin{corollary}\label{ladm}
A sequence $s_{\lambda_1}, \dots, s_{\lambda_{2t}}$
is $l$-admissible if and only if 
the tuple $\sigma=(\lambda_1;\dots;\lambda_{2t})$
is $l$-admissible. In that case
$s_{\lambda_{2t}}\dots s_{\lambda_1}(l)=l^{\sigma}$.
\end{corollary}

\begin{proposition}\label{peq} An integer $m\in S(l)$ if and only if
$m-l=2k>0$ 
and $k\subset_p m$. 
\end{proposition}

\begin{proof} Let $\sigma=(\lambda_1;\dots;\lambda_{2t})$
be an $l$-admissible tuple and 
$m=s_{\lambda_{2t}}\dots s_{\lambda_1}(l)$. 
Set $m^0=l$, 
$m^j=s_{\lambda_{2j}}\dots s_{\lambda_1}(l)$, and 
$k^j=(m^j-m^{j-1})/2$, $1\le j\le t$. 
Using Proposition~\ref{l29}, one deduces that $k^j\subset_p m^j$;
$m^j=l^{\sigma^j}$ with $\sigma^j=(\lambda_1;\dots;\lambda_{2j})$;
$k^j_i=0$ and $m^j_i=m^{j-1}_i$ for 
$i\notin[\lambda_{2j},\lambda_{2j-1}]$.
Therefore $k=k^1+\ldots+k^t\subset_p m$.

Assume now that 
$m-l=2k>0$ and $k\subset_p m$. Choose integers 
$\tau_1,\tau_2, \ldots, \tau_{2t}$ as follows. Set
$\tau_0=-1$. Assume that $\tau_{2j}$ is chosen. If there is no 
$i>\tau_{2j}$ such that $k_i=m_i\neq0$, we set $t=j$ and stop the 
process. Otherwise we choose for $\tau_{2j+1}$ minimal $i>\tau_{2j}$
with $k_i=m_i\neq0$. As $m>2k$, there exists $f>\tau_{2j+1}$ with 
$k_f\neq m_f$ (observe that in this case $k_f=0$). We choose 
minimal such $f$ for $\tau_{2j+2}$. 
Set $\lambda_{q}=\tau_{2t+1-q}$.
Since $k\subset_p m$, using Corollary~\ref{ladm}
and analyzing the $p$-adic expansions of $k$ and $m$,
one can conclude that the tuple $\sigma=(\lambda_1;\dots;\lambda_{2t})$
is $l$-admissible and $m=l^{\sigma}$.
\end{proof}

\begin{corollary}\label{ceq} Theorems~\ref{Ad1n} and \ref{dec2}
are equivalent, so 
Theorems~\ref{dec} and \ref{dec2} are valid in
characteristic 2 as well. 
\end{corollary}

Now we rewrite Adamovich's results \cite{A2} on the submodule structure
of the Weyl modules in our terms. 
We fix  $n$ and write $V_l$ and $\pi_m$ instead of 
$V^n_l$ and $\pi^n_m$. For $m\in S(l)$ or $m=l$ we denote
by $P_l(m)$ the smallest submodule of $V_l$ that has a composition 
factor $\pi_m$. Since $V_l$ is multiplicity-free, $P_l(m)$ is correctly 
defined and each submodule of $V_l$ is a sum of $P_l(m)$ for 
some $m$. Hence the submodule structure of $V_l$ is determined by the
inclusion relations between the submodules $P_l(m)$  (see also comments 
at the beginning of \cite{A2}). 
We shall write $\pi_m\prec\pi_q$ if $P_l(m)\subset P_l(q)$.
Let $\sigma=(\lambda_1;\dots;\lambda_{2t})$ be an
$l$-admissible tuple. For 
$m=l^{\sigma}$ 
set 
$$
\gQ_l(m)=\gQ_l(\sigma)=\bigcup_{j=1}^t[\lambda_{2j},\lambda_{2j-1}-1].
$$
Note that $\gQ_l(m)=\bigcup_{j=1}^t[\tau_{2j-1},\tau_{2j}-1]$
where
$\tau_i$ are
as in the proof of
Proposition~\ref{peq}. 
For instance, for $p=3$,
\begin{eqnarray*}
m&=&(0,\ 1,\ \underline{2,\ 2,\ 0},\ 1,\ 0,\ \underline{2,\ 1,\ 0,\ 0},\ 2,\
1), \\
k&=&(0,\ 0,\ \underline{2,\ 2,\ 0},\ 0,\ 0,\ \underline{2,\ 1,\ 0,\ 0},\ 0,\
0),
\end{eqnarray*}
we have $\gQ_{m-2k}(m)=[2,4]\cup[7,10]$. Put also $\gQ_l(l)=\emptyset$. 

\begin{theorem}\cite{A2} \label{Ad2b}  For $1\le m,q\le n+1$ and $m,q\in
S(l)\cup\{l\}$ 
the module 
$\pi_m\prec\pi_q$ (as composition factors of $V_l$) if and only if 
$\gQ_l(q)\subset\gQ_l(m)$.
\end{theorem}

\begin{remark}\label{re} Actually the sets $\gP_l(m)$ which are considered 
in \cite{A2} differ slightly  from $\gQ_l(m)$. For $m\in S(l)$ one has 
$
\gP_l(m)=\cup_{j=1}^t[\mu_{2j}+1,\mu_{2j-1}]
$ 
where $\mu_s=\lambda_s$ for $s<2t$,  $\mu_{2t}=\lambda_{2t}$ if 
$m\neq s_{\lambda_{2t-1}}\ldots s_{\lambda_1}(l)$, and  $\mu_{2t}=0$
otherwise. However, Lemma~\ref{l28} enables one to deduce that 
$\gP_l(m)\subset \gP_l(q)$ if and only if $\gQ_l(m)\subset \gQ_l(q)$. 
The crucial point is that $\lp(l)=\lp(m)$ for $m\in S(l)$.
\end{remark}

For $l$-admissible tuples
$\sigma=(\lambda_1;\dots;\lambda_{2t})$ 
and $\sigma'=(\lambda_1';\dots;\lambda_{2s}')$ 
we say that $\sigma\le\sigma'$ if 
there exists $f\le 2t,2s$ such that
$\lambda_i=\lambda_i'$ for $1\le i\le f$ and 
either $f=2t$, or $f<2t,2s$ and $\lambda_{f+1}<\lambda_{f+1}'$. 
It is convenient to assume that the empty tuple $\emptyset$
is $l$-admissible, $\emptyset\le\sigma$ for all $\sigma$,
$l^\emptyset=l$, and $\gQ_l(\emptyset)=\emptyset$.
The following is obvious.

\begin{lemma}\label{ll1}
Let $\sigma$ and $\sigma'$ be $l$-admissible tuples. 
Then $l^\sigma\le l^{\sigma'}$ if and only if
$\sigma\le\sigma'$.
\end{lemma}

Set $n'=n+1$. Construct an $l$-admissible tuple 
$\sigma^{\rm max}=(\mu_1;\dots;\mu_{2t})$ as follows.
Put $\mu_0=+\infty$. Assume that $\mu_{2j}$ is chosen.
Set $\mu=\mu_{2j}-1$.
If there is no $l$-admissible tuple $(\alpha;\beta)$ such that
$\mu\ge\alpha>\beta$ and 
$$
(l_0,\dots,l_{\mu})^{(\alpha;\beta)}
=(l_0,\dots,{\bar l}_{\beta}+1,{\bar l}_{\beta+1},
\dots,{\bar l}_{\alpha-1},l_{\alpha}+1,\dots,l_{\mu})
\le
(n_0',\dots,n_{\mu}'),
$$
we stop the process and set $t=j$ 
($\sigma^{\rm max}=\emptyset$ if $t=0$).
Otherwise we choose maximal such pair $(\alpha;\beta)$
(with respect to $\le$);
set $\mu_{2j+1}=\alpha$ and  $\mu_{2j+2}=\beta$;
and if 
$$
({\bar l}_{\beta}+1,{\bar l}_{\beta+1},\dots,{\bar l}_{\alpha-1},l_{\alpha}+1)
<(n_{\beta}',\dots,n_{\alpha}'),
$$
we stop the process and determine
$(\mu_{2j+3};\dots;\mu_{2t})$ as the maximal $l$-admissible
tuple with $\mu_{2j+3}<\beta$.
Obviously, $l^{\sigma^{\rm max}}$ is the maximal integer $m$
such that $\pi_m^n$ is a composition factor of $V_l^n$.

For $l$-admissible tuples $\sigma$ and $\sigma'$
we write $\sigma\prec\sigma'$ if and only if
$\gQ_l(\sigma)\supset\gQ_l(\sigma')$.
Using
Corollary~\ref{ladm}, Theorem~\ref{Ad2b}, and Lemma~\ref{ll1},
we get our main result on the structure of
fundamental Weyl modules.

\begin{theorem}\label{weylmain}
The map $\sigma\mapsto\pi_{l^\sigma}^n$ is a poset isomorphism
between the $l$-admissible tuples $\sigma\le\sigma^{\rm max}$ and the
composition factors of $V_l^n$ with the partial orders $\prec$.
\end{theorem}

If $l<n'$, we denote by $v$ the maximal integer
such that $l_v\ne n_v'$. If $l_v+1=n_v'$,
we denote by $u$ the maximal integer $<v$ such that
${\bar l}_u\ne n_u'$ setting $u=-1$ if 
${\bar l}_i=n_i'$ for all $i<v$. Put $s=\lp(l)$.

\begin{corollary}\label{irredweyl}
Let $n'=n+1$ and $1\le l\le n'$. 
Then $V_l^n$ is irreducible (i.e. $\sigma^{\rm max}=\emptyset$) 
if
and only if one of the following holds.

(1) $l=n'$;

(2) $l<n'$ and $s\ge v$;

(3) $l<n'$, $s<v$, $l_v+1=n_v'$, ${\bar l}_s\ge n_s'$;
$l_i=p-1$ and $n_i'=0$ for $s<i<v$.
\end{corollary}

\begin{proof} This follows from Proposition~\ref{l29} and 
Corollary~\ref{ladm}.
\end{proof}

\begin{corollary}\label{socleweyl}
Let $n'=n+1$ and $1\le l\le n'$.
The socle of $V_l^n$ is always simple. For reducible
$V_l^n$ it has the form 
$\pi_{l^\gamma}$ with $\gamma=(t;s)$ and $t$ as follows.

(1) $t=v$ if $s<v$ and either $l_v+1<n_v'$, or ${\bar l}_u<n_u'$;

(2) $t=w$ if $l_v+1=n_v'$; $u=-1$ or ${\bar l}_u>n_u'$; $s<w<v$; 
$l_w\ne p-1$; and $l_j=p-1$ for $w<j<v$.
\end{corollary}

\begin{proof} Applying Results~\ref{l29}, \ref{ladm},
and \ref{weylmain},
we conclude that $\pi_{l^\gamma}$ is a composition factor
of $V_l^n$ and for each $l$-admissible tuple $\tau\le\sigma^{\rm max}$
the set 
$\gQ_l(\tau)\subset[s,t-1]$, so $\pi_{l^{\gamma}}\prec\pi_{l^{\tau}}$.
\end{proof}


\section{Branching rules and the submodule structure of the restrictions}
\label{3}

In this section the main results of the article are proved. 
We shall need the following simple lemma.

\begin{lemma}\label{ldop1} Assume that $\dkl=1$ (i.e. $k\subset_p l+2k$).

(i) If $p^s\mid l+2k$, then $p^s\mid k$ and $p^s\mid l$.

(ii) If $p^s\mid l$, then $p^s\mid k$ and $p^s\mid l+2k$.
\end{lemma}

\begin{proof} One can assume that $s\ge 1$.

\noindent
$(i)$ Let $p^s\mid l+2k$. Since $k\subset_p l+2k$, we have $p^s\mid k$. This
implies that $p^s\mid l$.

\noindent
$(ii)$ Let $p^s\mid l$. Then $l_0=\dots=l_{s-1}=0$. 
Let $r=\lp(k)$. Assume that $r<s$. 
Since $k\subset_p l+2k$, we have $k_r=(2k)_r\ne 0$, which is impossible.
Therefore $r\ge s$, so $p^s\mid k$ and $p^s\mid l+2k$.
\end{proof}

As in Section~\ref{2} , we shall omit the superscript $n$ in our notation for
modules when it is known which
group is considered. Replacing $\om_i$ by $\pi_{n+1-i}$ and $W_i$ by
$V_{n+1-i}$, one immediately concludes that 
Theorem~\ref{br}(i) is equivalent to the following

\begin{theorem}\label{br2}
Let $1\le i\le n+1$ and $d=\lp(i)$. Then
$$
\pi_i^n\restr G_{n-1}\sim 
\pi_{i-1}+2\pi_i+\left(\sum_{t=0}^{d-1}2\pi_{i-1+2p^t}\right)+
\varepsilon\pi_{i-1+2p^d}
$$
where $\varepsilon=0$ if $i\equiv-p^d \pmod{p^{d+1}}$ and 
$\varepsilon=1$ otherwise.
\end{theorem}

\begin{proof} One can  rewrite the formula in Theorem~\ref{br2}
as follows.
\begin{equation}\label{f10}
\pi_i^n\restr G_{n-1}\sim 
\pi_{i-1}+2\pi_{i}+\sum_{t=0}^{\infty}b_t^i\pi_{i-1+2p^t}
\end{equation}
where
\begin{equation}\label{f101}
b_t^i=\left\{
\begin{array}{rl}
2, & i\equiv 0 \pmod{p^{t+1}}, \\
1, & i\equiv ap^t \pmod{p^{t+1}} \mbox{ and } 
a\not\equiv 0,-1 \pmod{p}, \\
0, & i\not\equiv 0 \pmod{p^t} \mbox{ or } i\equiv -p^t \pmod{p^{t+1}}.
\end{array}
\right.
\end{equation}
Recall that 
by convention $\pi_i^n=0$ for all $i>n+1$,
and $\pi_{n+1}^n$ is the trivial one-dimensional $G_n$-module. So 
(\ref{f10}) holds for $i\ge n+1$. 
Assume now that 
$1\le l < n+1$ and 
(\ref{f10}) is
valid for all $i>l$. We shall prove it for $i=l$. Then the theorem
will
follow by induction. 

It follows from \cite[Proposition~3.3.2 and Theorem~4.3.1]{D} that 
$V_l\restr G_{n-1}$ has a filtration by Weyl modules for $G_{n-1}$. 
Then the classical branching rules for characteristic 0 
\cite{Zhe}
and Theorem~\ref{dec2} imply
\begin{equation}\label{eqfi}
V_l^n\restr G_{n-1}\sim V_{l-1}+2V_l+V_{l+1}
\sim 
2{\cal V}+
\sum_{t=0}^{\infty}f_t^{l+2t-1}\pi_{l+2t-1}
\end{equation}
where ${\cal V}=\sum_{k=0}^{\infty}d_k^{l+2k}\pi_{l+2k}$,\quad
$f_0^{l-1}=d_0^{l-1}$, and 
$$
f_t^{l+2t-1}=d_t^{l+2t-1}+d_{t-1}^{l+2t-1}\quad\mbox{for }t\ge1.
$$
On the
other hand, by Theorem~\ref{dec2}, 
$
V_l\restr G_{n-1}\sim 
\sum_{k=0}^{\infty}d_k^{l+2k}(\pi_{l+2k}\restr G_{n-1}).
$
Since $d_0^l=1$ and the branching rules for $\pi_i$ with $i>l$ 
are assumed to satisfy (\ref{f10}), one can determine the 
branching of $\pi_l$.  
Therefore it suffices to check that the right part
of (\ref{eqfi}) is equal to 
$\sum_{k=0}^{\infty}d_k^{l+2k}U_{l+2k}$
where $U_i$ is the right part of (\ref{f10}).
The latter sum can be rewritten as follows: 
$$
2{\cal V}+
\sum_{k=0}^{\infty}d_k^{l+2k}(\pi_{l+2k-1}+
\sum_{s=0}^{\infty}b_s^{l+2k}\pi_{l+2k+2p^s-1}) 
=
2{\cal V}+
\sum_{t=0}^{\infty}e_t^{l+2t-1}\pi_{l+2t-1}
$$
where 
\begin{equation}\label{f13}
e_t^{l+2t-1}=d_t^{l+2t}+\sum_{k,s\ge0,\ k+p^s=t}d_k^{l+2k}b_s^{l+2k}.
\end{equation}
We have to show that $e_t^{l+2t-1}=f_t^{l+2t-1}$ for all $t\ge 0$. Note that
$f_t^{l+2t-1}\le2$.
We proceed by steps.

\medskip
\noindent
{\em Step 1.} {\em At most one summand in (\ref{f13}) is
nonzero. In particular, $e_t^{l+2t-1}\le2$}.

Assume that $d_{k}^{l+2k}b_{s}^{l+2k}\ne0$ and
$d_{k'}^{l+2k'}b_{s'}^{l+2k'}\ne0$ with $t=k+p^{s}=k'+p^{s'}$ and
$s>s'$. Since $b_{s}^{l+2k}\ne0$, we have $p^{s}\mid l+2k$. As
$d_{k}^{l+2k},d_{k'}^{l+2k'}\ne0$, 
by Lemma~\ref{ldop1}, $p^s$ divides $k$, $l$, and $k'$.
Hence $p^{s}\mid k+p^{s}-k'=p^{s'}$,
which yields a contradiction.

Now assume  that $d_t^{l+2t}\ne0$ and $d_k^{l+2k}b_s^{l+2k}\ne0$ with
$k+p^s=t$. As above, we get that $p^s$ divides $k$, $l$, and $t$.
Let $r=\lp(l)$. Then by
Lemma~\ref{ldop1}~(ii),  $p^r\mid k$ and $p^r\mid t$. Since $k+p^s=t$, we
have $r=s$, so $l_s\ne0$. Consider the following cases.

\medskip
\noindent
{\em Case 1.} $k\not\equiv0$ and $t\not\equiv0 \pmod{p^{s+1}}$. Then 
$l+2t\equiv t$ and 
$l+2k\equiv k \pmod{p^{s+1}}$, so $p^s=t-k\equiv2(t-k) \pmod{p^{s+1}}$,
which is impossible.

\medskip
\noindent
{\em Case 2.} $k\equiv0 \pmod{p^{s+1}}$. Then $t_s=1$. Since 
$l+2t\equiv t \pmod{p^{s+1}}$, we have $l_s=p-1$, so $(l+2k)_s=p-1$. This
implies that $b_s^{l+2k}=0$ and yields a contradiction.

\medskip
\noindent
{\em Case 3.} $t\equiv0 \pmod{p^{s+1}}$. Then $k_s=p-1$. Therefore
$(l+2k)_s=k_s=p-1$, so as above, $b_s^{l+2k}=0$. 

\medskip
\noindent
{\em Step 2.} {\em $e_t^{l+2t-1}=2$ if and only if $f_t^{l+2t-1}=2$
(equivalently, $d_t^{l+2t-1}=d_{t-1}^{l+2t-1}=1$)}. 

Assume that $e_t^{l+2t-1}=2$. By Step 1, this is equivalent to the following:
there exist $k,s\ge0$ with $k+p^s=t$ such that $b_s^{l+2k}=2$ and 
$d_k^{l+2k}=1$. Hence $p^{s+1}\mid l+2k$ and $k\subset_p l+2k$. By
Lemma~\ref{ldop1}~(i), $p^{s+1}$ divides $l$ and $k$. Note that
$l+2t-1=l+2k+p^s+(p^s-1)$. Therefore $t=k+p^s\subset_p l+2t-1$ and
$t-1=k+(p^s-1)\subset_p l+2t-1$, so $d_t^{l+2t-1}=d_{t-1}^{l+2t-1}=1$, as
required.

Assume now that $d_t^{l+2t-1}=d_{t-1}^{l+2t-1}=1$. Let $s=\lp(t)$. 
Then $t_s\ne0$ and $(t-1)_s=t_s-1$. Since
both $t$ and $t-1$ are contained in $l+2t-1$, we have $t_s=1$. Moreover, we
have $(l-1)_0=\dots=(l-1)_{s-1}=(t-1)_0=\dots=(t-1)_{s-1}=p-1$. 
Since $(l+2t-1)_s=t_s=1$, we get $(l-1)_s=p-1$. Hence $p^{s+1}\mid l$.
Set $k=t-p^s$. Then $p^{s+1}\mid k$, so $p^{s+1}\mid l+2k$ and
$b_s^{l+2k}=2$. It remains to observe that $k=t-p^s\subset_p
(l+2t-1)-p^s-(p^s-1)$,
so $d_k^{l+2k}=1$ and $e_t^{l+2t-1}=d_k^{l+2k}b_s^{l+2k}=2$.

\medskip
\noindent
{\em Step 3.} {\em If $e_t^{l+2t-1}\ne0$, then $f_t^{l+2t-1}\ne0$,
i.e. $t\subset_p l+2t-1$ or $t-1\subset_p l+2t-1$}.

Let $r=\lp(l)$. First assume  that
$d_t^{l+2t}=1$ (see (\ref{f13})), i.e. $t\subset_p l+2t$. Then by
Lemma~\ref{ldop1}~(ii), $p^r\mid t$. One easily checks that if 
$p^{r+1}\nmid t$, then $t-1\subset_p l+2t-1$, and if $p^{r+1}\mid t$, then 
$t\subset_p l+2t-1$, as required. 

Assume now that there exist $k,s\ge0$ with $k+p^s=t$ 
such that $d_k^{l+2k}=1$ and $b_s^{l+2k}\ne0$. 
By Lemma~\ref{ldop1}, $p^r\mid k$ and $\lp(l+2k)=r$. If $b_s^{l+2k}=2$, then
by Step 2, $f_t^{l+2t-1}=2\ne0$.
Hence we can assume that $b_s^{l+2k}=1$. By (\ref{f101}), $p^s\mid l+2k$
and $(l+2k)_s\ne0,p-1$, so $r=s$. Assume that $k_s=0$. Then
$(t-1)_s=0$ and $t-1=k+p^s-1\subset_p l+2k+p^s+(p^s-1)=l+2t-1$.
If $k_s\ne0$, we have $k_s=(l+2k)_s\ne p-1$, so
$t=k+p^s\subset_p l+2k+p^s+(p^s-1)=l+2t-1$, as required.

\medskip
\noindent
{\em Step 4.} {\em If $f_t^{l+2t-1}\ne0$, then $e_t^{l+2t-1}\ne0$}.

If $t=0$, then $e_0^{l-1}=1$, so assume that $t\ge1$. 
We have either $t\subset_p l+2t-1$, or $t-1\subset_p l+2t-1$. One needs to
show
that either $t\subset_p l+2t$ (i.e. $d_t^{l+2t}=1$), or 
there exist $k,s\ge0$ with $k+p^s=t$ 
such that $k\subset_p l+2k$, $p^s\mid l+2k$,
and $(l+2k)_s\ne p-1$ (i.e. $d_k^{l+2k}=1$ and $b_s^{l+2k}\ne0$). 
First assume that $t\subset_p l+2t-1$. Let $s=\lp(t)$. 
We have $(l+2t-1)_s=t_s\ne0$. If  $p^s\mid l+2t$,  then 
$k=t-p^s\subset_p l+2t-1-p^s-(p^s-1)=l+2k$,
$p^s\mid l+2t-2p^s=l+2k$, and $(l+2k)_s=(l+2t-1)_s-1\ne p-1$,
as required. 
If $p^s\nmid l+2t$, one gets $t\subset_p (l+2t-1)+1=l+2t$, as
desired.

Assume now that $t-1\subset_p l+2t-1$. Consider the following cases.

\medskip
\noindent
{\em Case 1.} $(t-1)_0=0$. If $(l+2t-1)_0=0$, then $t\subset_p l+2t$. Assume
that
$(l+2t-1)_0\ne0$. Set $s=0$, $k=t-1$. Then 
$k=t-1\subset_p (l+2t-1)-1=l+2k$ and $(l+2k)_0=(l+2t-1)_0-1\ne p-1$, as
required.

\medskip
\noindent
{\em Case 2.} $(t-1)_0\ne0,p-1$. Then $(l+2t-1)_0=(t-1)_0\ne p-1$, so
$t\subset_p l+2t$.

\medskip
\noindent
{\em Case 3.} $p\mid t$. Let $s=\lp(t)$.
Since $t-1\subset_p l+2t-1$, we have $p^s\mid l+2t$, so $p^s\mid l$. As 
$p^{s+1}\nmid t$, the integer $(t-1)_s\ne p-1$. If $(t-1)_s\ne0$, then 
$(l+2t-1)_s=(t-1)_s\ne p-1$, so $t\subset_p l+2t$. Assume now that
$(t-1)_s=0$.
This implies $t_s=1$. If $(l+2t-1)_s=0$, then $t\subset_p l+2t$. Therefore
one
can suppose that $(l+2t-1)_s\ge 1$. 
Then $k=t-p^s\subset_p (l+2t-1)-p^s-(p^s-1)=l+2k$,
$p^s\mid l+2t-2p^s=l+2k$, and $(l+2k)_s=(l+2t-1)_s-1\ne p-1$, as
required. 
\end{proof}

Now we investigate the submodule structure of the restriction
$\pi_i^n\restr G_{n-1}$. Let $n>1$ and $1\le i\le n$. 
As $\pi_i^n$ is the top composition factor of
$V_i^n$, it follows from (\ref{eqfi}) that $\pi_i^n\restr G_{n-1}$ is a 
quotient of the $G_{n-1}$-module
$V_i^n\restr G_{n-1}\sim V_{i-1}+2V_i+V_{i+1}$. Applying Smith's theorem
\cite{Sm} both to $V_i$ and $\pi_i$, we conclude that 
$V_i^n\restr G_{n-1}=V_i\oplus V_i\oplus V$ where
$V\sim V_{i-1}+V_{i+1}$, and 
$\pi_i^n\restr G_{n-1}=\pi_i\oplus\pi_i\oplus D$ where $D$ is a quotient
of $V$. 
Now Theorem~\ref{br}(ii) and Corollary~\ref{ccr} follow 
immediately from 

\begin{theorem}\label{soc2} Let $d=\lp i$. 
Set $\varepsilon=0$ if $i\equiv-p^d \pmod{p^{d+1}}$ and 
$\varepsilon=1$ otherwise;
$j_q=i-1+2p^q$.
Choose minimal $t\in\Ze^+$ 
such that $j_t>n$. Put $d'=\min\{d,t\}$. 
Then 
$$
\pi_{j_0}\prec_s\pi_{j_1}\prec_s\dots
\prec_s\pi_{j_{d'-1}}\prec_s
\pi_{i-1}\oplus \varepsilon\pi_{j_d}
\prec_s\pi_{j_{d'-1}}\prec_s
\dots\prec_s\pi_{j_1}\prec_s\pi_{j_0} 
$$
is the socle series of $D$. In particular, 
$D=\pi_{i-1}\oplus\varepsilon\pi_{j_d}$ if $d'=0$. 
\end{theorem}

\begin{proof} By Theorem~\ref{br2}, 
$
D\sim\pi_{i-1}+2\pi_{j_0}+\dots+2\pi_{j_{d'-1}}+
\varepsilon\pi_{j_d}.
$
It follows from Theorem~\ref{dec2} that the
factors $\pi_{i-1},\pi_{j_0},\dots,\pi_{j_{d'-1}}$ come
from $V_{i-1}$ and the factors 
$\pi_{j_0},\dots,\pi_{j_{d'-1}}$, and 
$\varepsilon\pi_{j_d}$ if nonzero 
come from $V_{i+1}$. Note that
$$
j_k=i-1+2p^k=i+p^k+(p^k-1)=i+(\underbrace{p-1,\dots,p-1}_k,1),
$$
so $\gQ^{i-1}(j_k)=[k,d-1]$ for all $0\le k\le d-1$.
Therefore by Theorem~\ref{Ad2b},
\begin{equation}\label{f33}
\pi_{j_0}\prec\pi_{j_1}\prec\dots
\prec\pi_{j_{d'-1}}\prec
\pi_{i-1}
\quad\mbox{in $V_{i-1}$}.
\end{equation}
Similarly, we get $\gQ^{i+1}(j_k)=[0,k-1]$ for
$1\le k<d$ (and for $k=d$ if $\varepsilon\ne0$ and $d>0$). Hence
\begin{equation}\label{f34}
\varepsilon\pi_{j_d}\prec\pi_{j_{d'-1}}\prec
\dots\prec\pi_{j_1}\prec\pi_{j_0}\quad\mbox{in $V_{i+1}$}.
\end{equation}
(Here the symbol $\prec$ is extended to the zero module in the 
natural way.)
Since
$\pi_i^n$ is selfdual, $D$ is selfdual also.
Let $D_1\prec_s\dots\prec_s D_m$ be the socle series
of $D$. 
Recall that $D$ has a filtration by quotients of $V_{i-1}$
and $V_{i+1}$.
As the factor $\pi_{i-1}$ has multiplicity 1 and
$D$ is selfdual,
(\ref{f33}) implies that $\pi_{i-1}$ is a factor of $D_q$ with 
$d'+1\le q\le m-d'$, so $m\ge 2d'+1$. 
If $\varepsilon\pi_{j_d}=0$,
then $m=2d'+1$
is the composition length
of $D$
and the theorem follows from (\ref{f33}) and (\ref{f34}).
Assume that $\varepsilon\pi_{j_d}\ne0$. As above, by 
the selfduality of $D$ and
(\ref{f34}),
$\pi_{j_d}$ is a factor of $D_{q'}$ with $d'+1\le q'\le m-d'$.
Assume that $q'\ne q$. Then $m=2d'+2$, so $D$ is uniserial,
which contradicts the selfduality of $D$.
Hence $q'=q$ and the theorem follows from (\ref{f33}) and (\ref{f34}).
\end{proof}

\begin{remark} 
Obviously, if 
$\varepsilon\pi_{j_d}=0$
(i.e. $d'<d$ or $i\equiv-p^d \pmod{p^{d+1}}$), then
the module $D$ is uniserial, so has exactly $2d'+2$ different
submodules. Since $D$ is selfdual,
one can easily observe that $D$ has exactly $2d+4$ different
submodules in the case where 
$\varepsilon\pi_{j_d}\ne0$
(i.e. 
$d=d'$ and $i\not\equiv-p^d \pmod{p^{d+1}}$).
\end{remark}


\section{Inductive systems and the transfer to an 
arbitrary field} 
\label{4}

In this section the inductive systems of fundamental representations
for $Sp_{\infty}(K)$ 
are classified and the results of the paper are transferred 
to an arbitrary field.  

\begin{proposition}\label{tr} 

{\rm(i)} Let $\Psi$ be an inductive system and
${\cal R}^{p^{t-1}}\subset\Psi$.
Then ${\cal R}^{p^t-1}\subset\Psi$.

{\rm(ii)} 
${\cal R}^{k}$ is an inductive system if and only if $k=p^t-1$, $t\ge1$.
\end{proposition}

\begin{proof} (i) Assume that ${\cal R}^{p^t-1}\not\subset\Psi$. Choose
maximal 
$l$ such that ${\cal R}^l\subset\Psi$. Then $p^{t-1}\le l<p^t-1$. 
Take minimal 
$s$ such that $l_s\neq p-1$ and set $i=p^s(l_s,l_{s+1}, \ldots)$. 
Then $i>0$ and ${\cal R}^i\subset\Psi$. One has $\lp(i)\ge s$. Moreover, 
if $\lp(i)=s$, then $i\not\equiv -p^s \pmod{p^{s+1}}$. Therefore 
Theorem~\ref{br2} implies that 
$\pi^{n-1}_{i-1+2p^s}$ is a composition factor of 
$\pi^n_i\restr G_{n-1}$ for $n\ge i-1+2p^s$,
and ${\cal R}^{i-1+2p^s}\subset\Psi$. As $i-1+2p^s>l$, we get 
a contradiction.

(ii) In view of (i), it suffices to verify that
${\cal R}^{p^t-1}$ is an inductive system. 
By Theorem~\ref{br2}, we need only to check that 
if $i<p^t-1$ and $s=\lp(i)$, then $i-1+2p^s\le p^t-1$ if
$i_s\neq p-1$, and $i-1+2p^{s-1}\le p^t-1$ if $i_s=p-1$ and $s>0$. But this 
is 
clear since $i\le (p-2)p^s+(p-1)p^{s+1}+\ldots+(p-1)p^{t-1}$ in the first
case and $i\le (p-1)p^s+(p-1)p^{s+1}+\ldots+(p-1)p^{t-1}$ in the second one. 
\end{proof} 

\begin{proof}[Proof of Theorem~\ref{ti}.] Theorem~\ref{br}(i) and
Proposition~\ref{tr} 
yield that ${\cal F}$, ${\cal L}^s$, ${\cal R}^{p^t-1}$, and 
${\cal L}^s\cup{\cal R}^{p^t-1}$, $s\ge0$, $t\ge1$, 
are inductive systems for 
$Sp_{\infty}(K)$. 

Let $\Psi=\{\Psi_i, \ i=1,2,\ldots\}$ be an inductive system of fundamental 
representations. It is clear that 
either for every $s,u\in\Ze^+$ there exist $n$ and $l$ such that 
$\om^n_l\in\Psi_n$, $l>s$, and $n+1-l>u$, or 
$\Psi\subset{\cal L}^s\cup{\cal R}^u$ for some $s$ and $u$. In the first case
we claim that $\Psi={\cal F}$. Indeed, fix $m$ and $l$, $0\le l\le m$. Then
one can choose $k$ and $n$ such that $\om^n_k\in\Psi_n$, $k\ge l$, and
$n-k\ge m-l$. Since $\Psi$ is an inductive system, Theorem~\ref{br}(i) implies
that $\om^{n+l-k}_l\in\Psi_{n+l-k}$ and $\om^m_l\in\Psi_m$. Hence 
$\Psi={\cal F}$.

Next, suppose that $\Psi\subset{\cal L}^s\cup{\cal R}^u$. Choose minimal $s$
and $u$ with this property assuming that $s=-1$ if  $\Psi\subset{\cal R}^u$
and $u=0$ if  $\Psi\subset{\cal L}^s$. (Observe that for all $s$ and $u$, 
$({\cal L}^s\cap{\cal R}^u)_n=\emptyset$ for $n$ large enough.) We
shall prove that  $\Psi={\cal L}^s\cup{\cal R}^u$ and $u=p^t-1$ with 
$t\in\Ze^+$ (in particular, $\Psi={\cal R}^u$ for $s=-1$ and 
$\Psi={\cal L}^s$ for $u=0$).  

First let $u>0$. We claim that ${\cal R}^u\subset\Psi$ and $u=p^t-1$. As 
$\Psi\not\subset{\cal L}^s$ and $\Psi$ and ${\cal L}^s$ 
are inductive systems, 
$\Psi_n\cap{\cal R}^u_n\ne\emptyset$ for infinitely many
integers
$n$. So there exists $v\le u$ such that $\pi^n_v\in\Psi_n$  for infinitely
many $n$. Choose maximal such $v$. Theorem~\ref{br2} yields that
$\pi^n_v\in\Psi_n$  for all $n\ge v-1$ and ${\cal R}^v\subset\Psi$. Now 
Proposition~\ref{tr} and the choice of $v$ imply that $v=p^t-1$ and 
${\cal R}^v$
is an inductive system. It remains to show that $v=u$. Suppose this is
not the case. As $\Psi\not\subset{\cal L}^s\cup{\cal R}^v$, there exist $l$
and $t$ such that $v<l\le u$, $t>s+l-1$, and $\pi^t_l\in\Psi_t$. Since
$\Psi$, ${\cal L}^s$, and 
${\cal R}^v$ are inductive systems, this implies
that for every $k>t$ there exists $\pi^k_{m_k}\in\Psi_k$
with $v<m_k\le u$ which contradicts the choice of $v$. Hence $v=u=p^t-1$ and 
${\cal R}^u\subset\Psi$. 

Now we show that ${\cal L}^s\subset\Psi$ if $s\ge0$. As 
$\Psi\not\subset{\cal L}^{s-1}\cup{\cal R}^u$, for some $n>s+u-1$ we have 
$\om^n_s\in\Psi_n$. Since $\Psi$, ${\cal R}^u$, and ${\cal L}^{s-1}$ for
$s\ge1$ are inductive systems, this forces  $\om^n_s\in\Psi_n$ for all 
$n\ge s$. Now Theorem~\ref{br} yields that ${\cal L}^s\subset\Psi$, as
desired.  
\end{proof}        

\begin{proposition}\label{lred} All theorems of the paper
hold for arbitrary $F\supset K$. 
\end{proposition}

\begin{proof} Since the restrictions of the fundamental representations of
a semisimple algebraic group over an algebraically closed field to relevant 
Chevalley groups over arbitrary subfields remain irreducible and can be
realized over these subfields, only
Theorems~\ref{br}(ii) (or \ref{soc2}),
\ref{Ad2b}, and \ref{weylmain} require some analysis. 
Let $M$ be the $G_n$-module $\om^{n+1}_i\restr G_n$
or $W_i^n$.

First assume that $F=\bar F$ is algebraically closed.
Set $H=Sp_{2n}(F)$.
Let $L$ be the Lie algebra of $H$. For a root $\alpha$ of $H$ and $t\in F$ 
denote by $x_{\alpha}(t)\in H$ and $X_{\alpha}\in L$ the root elements in $H$
and $L$ associated with $\alpha$.  It is well known that 
$x_{\alpha}(t)(m)=(1+tX_{\alpha})m$ for $m\in M$ and long $\alpha$ (see, for 
instance, \cite[Lemma 1]{PS}).  For 
$g\in G_n$ set $x^g_{\alpha}(t)=gx_{\alpha}(t)g^{-1}$ and  
$X^g_{\alpha}=gX_{\alpha}g^{-1}$. It is clear that $X^g_{\alpha}\in L$. 
It suffices to show that each $G_n$-submodule
$N\subset M$ is an $H$-submodule. Obviously, $x^g_{\alpha}(t)N=N$ 
and $X^g_{\alpha}N\subset N$ for all long roots $\alpha$, $g\in G_n$, and
$t\in K$. But this forces $x^g_{\alpha}(t)N=N$ for all $t\in F$.
However, using the commutator relations for the Chevalley groups of type $C$
(see, for instance, \cite[Lemma~15]{St}), one can deduce that the subgroup 
generated by all
$x^g_{\alpha}(t)$
with $g\in G_n$, $t\in F$, and long $\alpha$
coincides with $H$. (Here, in fact,  it suffices to make computations 
within subgroups of type $C_2$ and show that our subgroup contains all 
short root subgroups). Hence $N$ is an $H$-module, as desired. 

Now let $F\supset K$ be arbitrary. For
a finite dimensional $FG_n$-module $S$
set $\bar S=S\otimes_F\bar F$ and denote the socle
of $S$ by ${\rm soc}(S)$. Since
$
\dim{\rm Hom}_{FG_n}(E,S)=\dim{\rm Hom}_{{\bar F}G_n}(\bar E,\bar S)
$
for any $FG_n$-module $E$, we have
${\rm soc}(\bar S)=\overline{{\rm soc}(S)}$ if all composition factors
of $S$ are absolutely irreducible.
The same holds for other members of the socle series of $S$.
If $M=W_i^n$, then $M$ and $\bar M$ are multiplicity-free 
and their submodules are completely determined
by the sets of composition factors. Therefore
the arguments on socles allow us to conclude that
each submodule of $\bar M$ has the form
$\bar S$ for some submodule $S\subset M$.

Let $M=\om^{n+1}_i\restr G_n$ with $i,n\ge1$. 
Then the socle of $M$ contains a submodule 
$V\cong\om_{i-1}^n\oplus\om_{i-1}^n$.
Since $M$ and $V$ are selfdual and $M$ has only two composition
factors isomorphic to $\om_{i-1}^n$, there exists
a submodule $D$ of $M$ such that $M=V\oplus D$ and $\bar D$
is the unique submodule in $\bar M$ with
$\bar M=\bar V\oplus \bar D$. Now one can see that
the socle series of $D$ is determined by that of $\bar D$
and is described by Theorem~\ref{br}(ii).
\end{proof}


\end{document}